\documentclass[11pt]{article}
\usepackage[utf8]{inputenc}

\usepackage{textcomp}
\usepackage[russian, english]{babel}

\usepackage{amsmath, amsfonts, amssymb,amsthm}
\usepackage{verbatim}

\usepackage{float}
\usepackage{changepage}
\usepackage{array}
\usepackage{enumerate}
\usepackage{hyperref}
\usepackage{authblk}

\usepackage[nottoc,numbib]{tocbibind}

\usepackage{xcolor}

\hypersetup{colorlinks,linkcolor={blue},citecolor={red},urlcolor={blue}}

\usepackage{tikz}
\usetikzlibrary{calc, matrix, arrows,decorations.pathmorphing, positioning}
\numberwithin{equation}{section}
\providecommand{\keywords}[1]
{
  \small	
  \textbf{\textit{Keywords---}} #1
}

\theoremstyle{definition}

\newtheorem{theorem}{Theorem}[section]

\newtheorem{consequence}{Consequence}[section]

\newcommand{\beq}{\begin{equation}}
\newcommand{\ee}{\end{equation}}

\newcommand{\ben}{\begin{equation*}}
\newcommand{\een}{\end{equation*}}

\newcolumntype{M}[1]{>{\centering\arraybackslash}m{#1}}

\title{Arnold's trivium problem №8}

\author[]{Timur Kenzhaev \thanks{kenzhaev\_t\_d@mail.ru}}

\affil[]{Skolkovo Institute of Science and Technology, Moscow, Russia}

\date{}

\begin{document}

\pagenumbering{arabic}

\maketitle

\begin{abstract}
We give a full and detailed solution of eighth Arnold's trivium problem. We find critical points of a smooth function on a given one-parametric two-dimensional surface with Lagrange multipliers method. Basic Morse theory and Poincaré–Hopf theorem makes it possible to determine a genus of this surface and gives a beautiful training example of Morse surgery on a two-dimensional surface of genus $g = 6$.
\end{abstract}


\section{Introduction}
Great mathematician Vladimir Igorevich Arnold published his Mathematical Trivium \cite{Arnold_Trivium} in 1991 lamenting the decline in mathematical culture of Russian universities. It is a list of one hundred ``simple" mathematical problems constituting (according to Arnold) a theoretical minimum of physics student. While some of these problems are simple exercises, the others are highly nontrivial and cause interest of students around the world. Bright example is problem №2:
\\\\
\textbf{Problem 2.} Find the limit
\ben
\lim\limits_{x\to 0}\, \frac{\sin \tan x - \tan \sin x}{\arcsin\arctan x - \arctan\arcsin x}.
\een
Although this problem might be easy solved (especially using a computer) with Maclaurin expansion, original author's solution using analiticity seems to be exotic for a modern math student. Moreover, generalization of this problem on smooth functions fails \cite{Osipov_Limit}. Arnold's trivium has been translated into English (for example, \cite{Khesin_Tabachnikov}) and causes active discussions on English-speaking mathematical forums. There have been attempts to collect solutions and drafts of these problems \cite{Khesin_Tabachnikov}, \cite{MathOMan}, but there's still no systematized set of full solutions. This text deals with the following problem:
\\\\
\textbf{Problem 8.} How many maxima, minima, and saddle points does the function $x^4 + y^4 + z^4 + u^4 + v^4$ have on the surface $x + \ldots + v = 0,\, x^2 + \ldots + v^2 = 1,\, x_1^3 + \ldots + v^3 = C$? 
\\\\
The answer is this:
\begin{enumerate}
    \item $C \in \left(-\frac{1}{\sqrt{30}}, \frac{1}{\sqrt{30}} \right)$~---~20 minima, 30 maxima and 60 saddle points;
    \item $C = \pm \frac{1}{\sqrt{30}}$~---~the surface is singular, there are 10 singular points which are maxima; in smooth part there are 20 minima and 30 saddle points;
    \item $C\in \left(\frac{1}{\sqrt{30}}, \frac{3}{\sqrt{20}}\right)\cup\left(-\frac{3}{\sqrt{20}}, -\frac{1}{\sqrt{30}}\right)$~---~20 minima, 20 maxima, 30 saddle points;
    \item $C = \pm\frac{3}{\sqrt{20}}$~---~the surface degenerates into 5 points;
    \item $C\notin \left[-\frac{3}{\sqrt{20}}, \frac{3}{\sqrt{20}} \right]$~---~the surface is the empty set.
\end{enumerate}
We obtain this result using Lagrange multipliers method. With Poincaré–Hopf theorem and basic Morse theory we conclude that in case of $C \in \left(-\frac{1}{\sqrt{30}}, \frac{1}{\sqrt{30}} \right)$ the surface is the sphere with 6 handles, which splits into 5 spheres passing through critical values $C = \pm \frac{1}{\sqrt{30}}$. These spheres degenerate into points and disappear passing through critical values $C = \pm\frac{3}{\sqrt{20}}$.  
\\\\
\textbf{Acknowledgements}. The author is grateful to Prof. Maxim Kazaryan for valuable discussions and comments on the draft of this article.
\section{Parameter $C$ variation}
Denote surface
\ben
\begin{cases}
& x + y + z + u + v = 0,
\\
& x^2 + y^2 + z^2 + u^2 + v^2 = 1
\end{cases}
\een
as $S^3$ (it is smooth compact three-dimensional surface diffeomorphic to three-dimensional sphere). As soon as it's compact, function $p_3 = x^3 + y^3 + z^3 + u^3 + v^3$ is bounded on $S^3$. Denote $p_i = x^i + y^i + z^i + u^i + v^i$ and $\mathcal{L}_1 = p_3 - \lambda_2\,p_2 - \lambda_1\,p_1$. Using Lagrange multipliers method we obtain a system of equations:
\ben
\begin{cases}
&d\mathcal{L}_1 = 0,
\\
&x + y + z + u + v = 0,
\\
&x^2 + y^2 + z^2 + u^2 + v^2 = 1.
\end{cases}
\Longleftrightarrow
\begin{cases}
&3x^2 - 2\lambda_2 x - \lambda_1 = 0,
\\
&3y^2 - 2\lambda_2 y - \lambda_1 = 0,
\\
&3z^2 - 2\lambda_2 z - \lambda_1 = 0,
\\
&3u^2 - 2\lambda_2 u - \lambda_1 = 0,
\\
&3v^2 - 2\lambda_2 v - \lambda_1 = 0,
\\
&x + y + z + u + v = 0,
\\
&x^2 + y^2 + z^2 + u^2 + v^2 = 1.
\end{cases}
\een
Coordinates of solution must be solutions of a quadratic equation, so there are three types of solutions (up to a permutation of coordinates):
\begin{enumerate}
\item $x = y = z = u = v$, there is no such solutions;
\item $x = y = z = u \neq v$, then $x = \pm \frac{1}{\sqrt{20}}, v = \mp\frac{4}{\sqrt{20}}$, $C = \mp\frac{3}{\sqrt{20}} \approx \mp 0,67082$;
\item $x = y = z \neq u = v$, then $x = \pm\frac{2}{\sqrt{30}}, v = \mp\frac{3}{\sqrt{30}}$, $C = \mp\frac{1}{\sqrt{30}}\approx \mp 0,18257$.
\end{enumerate}
It's clear that values $\mp\frac{3}{\sqrt{20}}$ are minimal and maximal, therefore, system
\beq
\label{Surface_main}
\begin{cases}
& x + y + z + u + v = 0,
\\
& x^2 + y^2 + z^2 + u^2 + v^2 = 1,
\\
& x^3 + y^3 + z^3 + u^3 + v^3 = C,
\end{cases}
\ee
has no solutions in case of $C\notin \left[-\frac{3}{\sqrt{20}}, \frac{3}{\sqrt{20}} \right]$ and has $\frac{5!}{4! 1!} = 5$ solutions in case of $C = \mp\frac{3}{\sqrt{20}}$.
\section{Surface smoothness}
Denote surface \eqref{Surface_main} as $A_3^C$. We need to check whether $A_3^C \subset\mathbb{R}^5$ is smooth to apply Lagrange multipliers method. From the Vandermonde determinant considerations it's easy to see that point $(x, y, z, u ,v)\in A_3^C$ is singular iff coordinates of this point take less than three values. Possible cases are similar to the previous paragraph:
\begin{enumerate}
\item $x = y = z = u = v$, there are no such points;
\item $x = y = z = u \neq v$, then $x = \pm \frac{1}{\sqrt{20}}, v = \mp\frac{4}{\sqrt{20}}$, $C = \mp\frac{3}{\sqrt{20}}$;
\item $x = y = z \neq u = v$, then $x = \pm\frac{2}{\sqrt{30}}, v = \mp\frac{3}{\sqrt{30}}$, $C = \mp\frac{1}{\sqrt{30}}$.
\end{enumerate}
Thus, in case of $C\in \left(-\frac{3}{\sqrt{20}}, \frac{3}{\sqrt{20}}\right)\setminus \left\{\mp\frac{1}{\sqrt{30}}\right\}$ the surface is smooth, in case of $C = \mp\frac{1}{\sqrt{30}}$ it has  $\frac{5!}{3! 2!} = 10$ singular points, and in case of $C = \mp\frac{3}{\sqrt{20}}$ it degenerates into 5 points.
\section{Critical points}
We need to find critical points of $p_4$ on $A_3^C$ in the smooth case (or in smooth part of the surface in the singular cases) with Lagrange multipliers method. Denote $\mathcal{L}_2 = p_4 - \lambda_3\,p_3 - \lambda_2\,p_2 - \lambda_1\,p_1$. Corresponding system is
\ben
\label{System_Main}
\begin{cases}
&d\mathcal{L}_2 = 0,
\\
&x + y + z + u + v = 0,
\\
&x^2 + y^2 + z^2 + u^2 + v^2 = 1,
\\
&x^3 + y^3 + z^3 + u^3 + v^3 = C.
\end{cases}
\Longleftrightarrow
\begin{cases}
&4x^3 - 3\lambda_3 x^2 - 2\lambda_2 x - \lambda_1 = 0,
\\
&4y^3 - 3\lambda_3 y^2 - 2\lambda_2 y - \lambda_1 = 0,
\\
&4z^3 - 3\lambda_3 z^2 - 2\lambda_2 z - \lambda_1 = 0,
\\
&4u^3 - 3\lambda_3 u^2 - 2\lambda_2 u - \lambda_1 = 0,
\\
&4v^3 - 3\lambda_3 v^2 - 2\lambda_2 v - \lambda_1 = 0,
\\
&x + y + z + u + v = 0,
\\
&x^2 + y^2 + z^2 + u^2 + v^2 = 1,
\\
&x^3 + y^3 + z^3 + u^3 + v^3 = C.
\end{cases}
\een
Coordinates $x, y, z, u, v$ of system \eqref{System_Main} solution $(x, y, z, u, v, \lambda_1, \lambda_2, \lambda_3)$  are roots of a cubic equation \\${4t^3 - 3\lambda_3 t^2 - 2\lambda_2 t - \lambda_1 = 0}$. It means that $(x_0, y_0, z_0, u_0, v_0, \lambda_1, \lambda_2, \lambda_3)$ is a solution iff coordinates of point $(x_0, y_0, z_0, u_0, v_0)\in A^{C}_3$ take less than four values. Indeed, otherwise $p_1, p_2, p_3, p_4$ can be complemented to curvilinear coordinate system in some neighbourhood of this point, thus, it's not a critical point. Lagrange multipliers then can be calculated via Viet theorem.
\\
We have already mentioned that there is no points of $A_3^C$ with $x = y = z = u = v$. If coordinates take two values it's either singular points in case of $C = \pm\frac{1}{\sqrt{30}}$ or isolated points in case of $C =\pm\frac{3}{\sqrt{20}}$. Then up to permutations, there are two types of critical points:
\begin{enumerate}
    \item $x = y\neq z = u\neq v$;
    \item $x = y = z\neq u\neq v$.
\end{enumerate}
\subsection{Points of the first type}
Denote $x = y = a \neq z = u = b \neq v = c$. We need to solve a system:
\ben
\begin{cases}
& 2a + 2b + c = 0;
\\
& 2a^2 + 2b^2 + c^2 = 1;
\\
& 2a^3 + 2b^3 + c^3 = C.
\end{cases}
\een
Expressing $c = -2(a + b)$ from the first equation we obtain a system of two equations
\ben
\begin{cases}
& 2a^2 + 2b^2 + 4(a + b)^2 = 1;
\\
& 2a^3 + 2b^3 - 8(a + b)^3 = C.
\end{cases}
\een
Solutions of this system are inetrsection points of conic and cubic curves. With natural change of variables $a = \frac{\tilde{a} + \tilde{b}}{2}, b = \frac{\tilde{a} - \tilde{b}}{2}$ the system is reduced to
\ben
\begin{cases}
&5\tilde{a}^2 + \tilde{b}^2 = 1,
\\
&\frac{3}{2}\, \tilde{a}\,\tilde{b}^2 - \frac{15}{2}\, \tilde{a}^3 = C.
\end{cases}
\een
Expressing $\tilde{b}$ from the first equation we obtain $a = \frac{t \pm \sqrt{1 - 5t^2}}{2}, b = \frac{t \mp \sqrt{1 - 5t^2}}{2}, c~=~-2t$, where $t$ is root of the cubic equation $15t^3 - \frac{3}{2}\, t + C = 0$, s.t. $t\in \left[-\frac{1}{\sqrt{5}}, \frac{1}{\sqrt{5}}\right]$. 
\\
Before exploring the number of solutions, we need to understand which their properties will be used. The next step is to calculate restriction of $d^2\mathcal{L}_2$ on a tangent space to the critical points. Tangent space is given by the system of equations:
\beq
\label{Tangent_First_System}
\begin{cases}
&dx + dy + dz + du + dv = 0,
\\
&a(dx + dy) + b(dz + du) + c\,dv = 0,
\\
&a^2(dx + dy) + b^2(dz + du) + c^2\, dv = 0.
\end{cases}
\ee
As soon as numbers $a, b, c$ are pairwise distinct, system \eqref{Tangent_First_System} is equivalent to $dx + dy = dz + du = dv = 0$. Then $d^2\mathcal{L}_2 = 2\,P_3'(a)dx^2 + 2\,P_3'(b)dz^2$, where $P_3(x) = 4x^3 - 3\lambda_3\,x^2 - 2\lambda_2\,x - \lambda_1 = 0$. Thus, we need to know the signs of $P_3$ derivative in its roots to determine the definition of $d^2\mathcal{L}_2$. As soon as $P_3$ has three distinct real roots $x_1 < x_2 < x_3$ and positive leading coefficient, $P'(x_1) > 0, P'(x_2) < 0, P'(x_3) > 0$. It follows immediately that critical points of the first type are nondegenerate. Moreover, they are either saddle points or maxima.
\\
We need to explore the relative positions of $a = \frac{t \pm \sqrt{1 - 5t^2}}{2}, b = \frac{t \mp \sqrt{1 - 5t^2}}{2}, c~=~-2t$, where $t\in \left[-\frac{1}{\sqrt{5}}, \frac{1}{\sqrt{5}}\right]$. Values $t = \pm \frac{1}{\sqrt{5}}, \pm \frac{1}{\sqrt{30}}$ do not suit us as soon as $a, b$ and $c$ fail to be pairwise distinct. Solving of inequalities gives:
\begin{enumerate}
\item with $t\in \left(-\frac{1}{\sqrt{5}}, -\frac{1}{\sqrt{30}}\right)$ relative positions of the roots are $a, b, c $ (or $b, a, c$), every $t$ within these boundaries gives $30$ saddle points; 
\item with $t\in \left(-\frac{1}{\sqrt{30}}, \frac{1}{\sqrt{30}}\right)$ relative positions of the roots are $a, c, b $ (or $b, c, a$), every $t$ within these boundaries gives $30$ maxima;
\item with $t\in \left(\frac{1}{\sqrt{30}}, \frac{1}{\sqrt{5}}\right)$ relative positions of the roots are $c, a, b $ (or $c, b, a$), every $t$ within these boundaries gives $30$ saddle points.
\end{enumerate}
Now we need to investigate number of polynomial $15t^3 - \frac{3}{2}\, t + C$ roots within intervals $\left(-\frac{1}{\sqrt{5}}, -\frac{1}{\sqrt{30}}\right)$, $\left(-\frac{1}{\sqrt{30}}, \frac{1}{\sqrt{30}}\right)$, $\left(\frac{1}{\sqrt{30}}, \frac{1}{\sqrt{5}}\right)$ while $C\in \left(-\frac{3}{\sqrt{20}}, \frac{3}{\sqrt{20}}\right)$. Simple graphical considerations give

\begin{itemize}
\item $C\in\left(-\frac{1}{\sqrt{30}}, \frac{1}{\sqrt{30}}\right)$~---~1 root in each of the three intervals;
\item $C = \frac{1}{\sqrt{30}}$~---~1 root in $\left(-\frac{1}{\sqrt{5}}, -\frac{1}{\sqrt{30}}\right)$, 1 root equal to $\frac{1}{\sqrt{30}}$;
\item $C\in\left(\frac{1}{\sqrt{30}}, \frac{3}{\sqrt{20}}\right)$~---~1 root in $\left(-\frac{1}{\sqrt{5}}, -\frac{1}{\sqrt{30}} \right)$;
\item $C = -\frac{1}{\sqrt{30}}$~---~1 root in $\left(\frac{1}{\sqrt{30}}, \frac{1}{\sqrt{5}}\right)$, 1 root equal to $-\frac{1}{\sqrt{30}}$;
\item $C\in\left(-\frac{3}{\sqrt{20}}, -\frac{1}{\sqrt{30}}\right)$~---~1 root in $\left(\frac{1}{\sqrt{30}} ,\frac{1}{\sqrt{5}}\right)$.
\end{itemize}

\subsection{Points of the second type}
Denote $x = y = z = \alpha \neq u = \beta \neq v = \gamma$. We need to solve a system:
\ben
\begin{cases}
& 3\alpha + \beta + \gamma = 0,
\\
& 3\alpha^2 + \beta^2 + \gamma^2 = 1,
\\
& 3\alpha^3 + \beta^3 + \gamma^3 = C.
\end{cases}
\een
Expressing $\alpha$ from the first equation we obtain a system
\ben
\begin{cases}
& \frac{1}{3}(\beta + \gamma)^2 + \beta^2 + \gamma^2 = 1,
\\
& -\frac{1}{9}(\beta + \gamma)^3 + \beta^3 + \gamma^3 = C.
\end{cases}
\een
With natural change of variables $\beta = \frac{\tilde{\beta} + \tilde{\gamma}}{2}, \gamma = \frac{\tilde{\beta} - \tilde{\gamma}}{2}$ the system is reduced to 
\ben
\begin{cases}
&5\tilde{\beta}^2 + 3\tilde{\gamma}^2 = 6 ,
\\
&\frac{5}{36} \tilde{\beta}^3 + \frac{3}{4} \tilde{\beta}\,\tilde{\gamma}^2 = C.
\end{cases}
\een
Similar to the previous paragraph we obtain $\alpha = -\frac{1}{3}\, t, \beta = \frac{t \pm \sqrt{2 - \frac{5}{3}\,t^2}}{2}, \gamma = \frac{t \mp \sqrt{2 - \frac{5}{3}\,t^2}}{2}$, where $t$ is a root of the cubic equaiton $\frac{10}{9}\,t^3 - \frac{3}{2}\,t + C = 0$, s.t. $t\in \left[-\sqrt{\frac{6}{5}}, \sqrt{\frac{6}{5}}\,\right]$.
\\
Tangent space is given by
\beq
\label{Tangent_Second_System}
\begin{cases}
&dx + dy + dz + du + dv = 0,
\\
&\alpha(dx + dy + dz) + \beta\,du + \gamma\,dv = 0,
\\
&\alpha^2(dx + dy + dz) + \beta^2\,du + \gamma^2\,dv = 0.
\end{cases}
\ee
As soon as numbers $\alpha, \beta$ and $\gamma$ are pairwise distinct, system \eqref{Tangent_Second_System} is equivalent to $dx + dy + dz = du = dv = 0$, the second differential on tangent space is $P_3'(\alpha)(dx^2 + dy^2 + dz^2)$ with restriction on $dx + dy + dz = 0$. Thus, critical points of the second type are nondegenerate and can be either maxima or minima. We need to explore the relative positions of $\alpha = -\frac{1}{3}\, t, \beta = \frac{t \pm \sqrt{2 - \frac{5}{3}\,t^2}}{2}, \gamma = \frac{t \mp \sqrt{2 - \frac{5}{3}\,t^2}}{2}$ where $t\in \left[-\sqrt{\frac{6}{5}}, \sqrt{\frac{6}{5}}\,\right]$. Values $t = \pm \sqrt{\frac{6}{5}}, \pm \frac{3}{2\sqrt{5}}$ do not suit us as soon as $\alpha, \beta, \gamma$ fail to be distinct. Solving of inequalities gives:

\begin{enumerate}
\item with $t\in \left(-\sqrt{\frac{6}{5}}, -\frac{3}{2\sqrt{5}}\right)$ relative positions of the roots are $\gamma, \beta, \alpha $ (or $\beta, \gamma, \alpha$), every t within these boundaries gives $20$ maxima; 
\item with $t\in \left(-\frac{3}{2\sqrt{5}}, \frac{3}{2\sqrt{5}}\right)$ relative positions of the roots are $\beta, \alpha, \gamma $ (or $\gamma, \alpha, \beta$), every t within these boundaries gives $20$ minima;
\item with $t\in \left(\frac{3}{2\sqrt{5}}, \sqrt{\frac{6}{5}}\right)$ relative positions of the roots are $\alpha, \beta, \gamma $ (or $\alpha, \gamma, \beta$), every t within these boundaries gives $20$ maxima.
\end{enumerate}
Now we need to investigate number of polynomial $\frac{10}{9}\,t^3 - \frac{3}{2}\, t + C$ roots within intervals $\left(-\sqrt{\frac{6}{5}}, -\frac{3}{2\sqrt{5}}\right)$, $\left(-\frac{3}{2\sqrt{5}}, \frac{3}{2\sqrt{5}}\right)$, $\left(\frac{3}{2\sqrt{5}}, \sqrt{\frac{6}{5}}\right)$ while $C\in C\in \left(-\frac{3}{\sqrt{20}}, \frac{3}{\sqrt{20}}\right)$. Simple graphical considerations give

\begin{itemize}
\item $C\in\left(-\frac{1}{\sqrt{30}}, \frac{1}{\sqrt{30}}\right)$~---~1 root in $\left(-\frac{3}{2\sqrt{5}}, \frac{3}{2\sqrt{5}}\right)$;
\item $C = \frac{1}{\sqrt{30}}$~---~1 root in $\left(-\frac{3}{2\sqrt{5}}, \frac{3}{2\sqrt{5}}\right)$, 1 root equal to $\sqrt{\frac{6}{5}}$;
\item $C\in\left(\frac{1}{\sqrt{30}}, \frac{3}{\sqrt{20}}\right)$~---~1 root in $\left(-\frac{3}{2\sqrt{5}}, \frac{3}{2\sqrt{5}}\right)$, 1 root equal to $\left(\frac{3}{2\sqrt{5}}, \sqrt{\frac{6}{5}}\right)$;
\item $C = -\frac{1}{\sqrt{30}}$~---~1 root in $\left(-\frac{3}{2\sqrt{5}}, \frac{3}{2\sqrt{5}}\right)$, 1 root equal to $-\sqrt{\frac{6}{5}}$;
\item $C\in\left(-\frac{3}{\sqrt{20}}, -\frac{1}{\sqrt{30}}\right)$~---~1 root in $\left(-\frac{3}{2\sqrt{5}}, \frac{3}{2\sqrt{5}}\right)$, 1 root in $\left(-\sqrt{\frac{6}{5}}, -\frac{3}{2\sqrt{5}}\right)$;
\end{itemize}

\subsection{Non-smooth case}
Above results show that in case of $C = \pm\frac{1}{\sqrt{30}}$ there is no maxima in smooth part of $A^C_{3}$, then from symmetry and compactness considerations all $10$ singular points are maxima. 
\subsection{Results and consequences}
Summarizing all of the above, we obtain the answer:
\begin{enumerate}
    \item $C \in \left(-\frac{1}{\sqrt{30}}, \frac{1}{\sqrt{30}} \right)$~---~20 minima, 30 maxima and 60 saddle points;
    \item $C = \pm \frac{1}{\sqrt{30}}$~---~the surface is singular, there are 10 singular points which are maxima; in smooth part there are 20 minima and 30 saddle points;
    \item $C\in \left(\frac{1}{\sqrt{30}}, \frac{3}{\sqrt{20}}\right)\cup\left(-\frac{3}{\sqrt{20}}, -\frac{1}{\sqrt{30}}\right)$~---~20 minima, 20 maxima, 30 saddle points;
    \item $C = \pm\frac{3}{\sqrt{20}}$~---~the surface degenerates into 5 points;
    \item $C\notin \left[-\frac{3}{\sqrt{20}}, \frac{3}{\sqrt{20}} \right]$~---~the surface is the empty set.
\end{enumerate}
\subsection{Genus of the surface}
In case of $C = 0$ one can notice $A^0_3$ is a connected surface. Indeed, function $p_4$ has 20 minima on $A^0_3$ with coordinates
\ben
x = y = z := \alpha \neq u := \beta \neq v := \gamma \:\text{ with }\: \alpha = 0,\, \beta = \pm \frac{1}{\sqrt{2}},\, \gamma = \mp \frac{1}{\sqrt{2}} \quad +\text{permutations}.
\een
Consider curve $ \gamma\colon [-1, 1] \to A^0_3$ with
\ben
\gamma(t) = 
\begin{pmatrix}
\sqrt{\frac{1 - t^2}{2}}, &-\sqrt{\frac{1 - t^2}{2}}, &0, &\frac{t}{\sqrt{2}}, & -\frac{t}{\sqrt{2}}
\end{pmatrix}
^T
\een
This curve connects points $\left(0, 0, 0, \frac{1}{\sqrt{2}}, -\frac{1}{\sqrt{2}}\right),
\left(\frac{1}{\sqrt{2}}, -\frac{1}{\sqrt{2}}, 0, 0, 0\right) \text{ and } \left(0, 0, 0, \frac{1}{\sqrt{2}}, -\frac{1}{\sqrt{2}}\right)$. It's easy to see that using such curves and permutations of coordinates one can connect any two minima of $p_4$ on $A_3^0$, which means $A_3^0$ is connected.
Then $\chi(A^0_3) = 20 - 60 + 30 = -10$ and $g(A^0_3) = 6$.
To define the type of the surface for $C\neq 0$ we need some basic facts from Morse theory.
For a smooth manifold $M$ and function $f\colon M\to\mathbb{R}$ denote set $\{p\in M \:|\: f(p) \leq a\}$ as $M^a$. 
\begin{theorem}
\cite{Morse_Theory} Let $f$ be a smooth real-valued function on a manifold $M$. Let $a < b$ and suppose that the set $f^{-1}\left([a, b]\right)$, consisting of all $p\in M$ with $a\leq f(p)\leq b$, is compact, and contains no critical points of $f$. Then $M^{a}$ is diffeomorphic to $M^b$. Furthermore, $M^a$ is a deformation retract of $M^b$, so that $M^a\to M^b$ is a homotopy equivalence.
\end{theorem}
\begin{consequence}
In case of $C \in \left(-\frac{1}{\sqrt{30}}, \frac{1}{\sqrt{30}} \right)$ $A_3^C$ is two-dimensional smooth compact surface of genus $6$.
\end{consequence}
\begin{theorem}
\cite{Morse_Theory} Let $f\colon M\to\mathbb{R}$ be a smooth function, and let $p$ be a non-degenerate critical point with index $\lambda$. Setting $f(p) = c$, suppose that $f^{-1}\left([c - \varepsilon, c + \varepsilon]\right)$ is compact, and contains no critical point of $f$ other than $p$, for some $\varepsilon > 0$. Then, for all sufficiently small $\varepsilon$, the set $M^{c + \varepsilon}$ has the homotopy type of $M^{c - \varepsilon}$ with a $\lambda$-cell attached.
\end{theorem}
\begin{consequence}
In case of $C\in \left(\frac{1}{\sqrt{30}}, \frac{3}{\sqrt{20}}\right)\cup\left(-\frac{3}{\sqrt{20}}, -\frac{1}{\sqrt{30}}\right)$ $A_3^{C}$ is the disjoint union of $5$ two-dimensional spheres. In case of $C = \pm \frac{1}{\sqrt{30}}$ $A_3^C$ is $5$ two-dimensional spheres which mutually intersect at one point having $C_5^2 = 10$ singular points. 
\end{consequence}

\bibliography{bibliography}{}

\begin{thebibliography}{MSW69}

\bibitem[Arn91]{Arnold_Trivium}
V.~I. Arnol'd.
\newblock {A mathematical trivium}.
\newblock {\em Russian Math. Surveys}, 46(1):~271--278, 1991.

\bibitem[KT14]{Khesin_Tabachnikov}
Boris Khesin and Serge Tabachnikov.
\newblock {\em {ARNOLD: Swimming Against the Tide}}.
\newblock Amer. Math. Soc., 2014.

\bibitem[Mat]{MathOMan}
\url{https://www.mathoman.com/1611-la-collection-d-exercices-de-vladimir-arnold/#co}.

\bibitem[MSW69]{Morse_Theory}
J.~Milnor, M.~SPIVAK, and R.~WELLS.
\newblock {\em Morse Theory. (AM-51), Volume 51}.
\newblock Princeton University Press, 1969.

\bibitem[Osi18]{Osipov_Limit}
N.~N. Osipov.
\newblock {{\foreignlanguage{russian}{Знаменитый предел Арнольда}}}.
\newblock {\em Mat. Pros., Ser.}, 3(22):~211--215, 2018.

\end{thebibliography}

\end{document}